\documentclass[11pt]{article}
\usepackage{amsfonts}
\usepackage{amssymb}
\usepackage{amsmath}
\def\sqr#1#2{{\vcenter{\vbox{\hrule height.#2pt
              \hbox{\vrule width.#2pt height#1pt \kern#1pt \vrule
width.#2pt}
              \hrule height.#2pt}}}}
\def\signed #1{{\unskip\nobreak\hfil\penalty50
              \hskip2em\hbox{}\nobreak\hfil#1
              \parfillskip=0pt \finalhyphendemerits=0 \par}}
\def\endpf{\signed {$\sqr69$}}
\def\dbE{{\mathbb{E}}}
\def\dbF{{\mathbb{F}}}

\def\dbM{{\mathbb{M}}}
\def\dbN{{\mathbb{N}}}

\def\dbP{{\mathbb{P}}}

\def\dbR{{\mathbb{R}}}

\def\dbW{{\mathbb{W}}}

\def\d{\delta}
\def\e{\varepsilon}

\def\t{\times}

\def\o{\omega}

\def\3n{\negthinspace \negthinspace \negthinspace }
\def\2n{\negthinspace \negthinspace }
\def\1n{\negthinspace }
\def\ns{\noalign{\smallskip} }
\def\ns{\noalign{\medskip} }
\def\ds{\displaystyle}
\def\O{\Omega}

\def\cF{{\cal F}}

\def\cM{{\cal M}}

\def\cl{{\cal l}}
\def\no{\noindent}

\def\ms{\medskip}
\def\bs{\bigskip}
\def\q{\quad}
\def\qq{\qquad}
\def\hb{\hbox}

\def\lan{\mathop{\langle}}
\def\ran{\mathop{\rangle}}

\def\wt{\widetilde}
\def\cd{\cdot}

\def\ae{\hbox{\rm a.e.{ }}}
\def\as{\hbox{\rm a.s.{ }}}

\def\cl{\overline}

\def\({\Big (}
\def\){\Big )}
\def\[{\Big[}
\def\]{\Big]}

\def\={\buildrel \triangle \over =}

\def\be{\begin{equation}}
\def\bel{\begin{equation}\label}
\def\ee{\end{equation}}
\def\bea{\begin{eqnarray}}
\def\eea{\end{eqnarray}}
\def\bt{\begin{theorem}}
\def\et{\end{theorem}}
\def\bc{\begin{corollary}}
\def\ec{\end{corollary}}
\def\bl{\begin{lemma}}
\def\el{\end{lemma}}
\def\bp{\begin{proposition}}
\def\ep{\end{proposition}}
\def\br{\begin{remark}}
\def\er{\end{remark}}
\def\ba{\begin{array}}
\def\ea{\end{array}}
\def\bd{\begin{definition}}
\def\ed{\end{definition}}

\newtheorem{lemma}{Lemma}[section]
\newtheorem{remark}{Remark}[section]

\newtheorem{theorem}{Theorem}[section]
\newtheorem{corollary}{Corollary}[section]

\newtheorem{definition}{Definition}[section]
\newtheorem{proposition}{Proposition}[section]

\makeatletter
   
   \@addtoreset{equation}{section}
\makeatother

\begin{document}

\title{\bf Well-posedness of Backward Stochastic Differential Equations with General Filtration}
\author{Qi L\"{u}\thanks{School of Mathematical Sciences, University of Electronic Science and Technology of China, Chengdu 610054, China;
and School of Mathematics, Sichuan University, Chengdu 610064,
China. {\small\it e-mail:} {\small\tt luqi59@163.com}. \ms}
 ~~~
  and~~~
Xu Zhang\thanks{Key Laboratory of Systems and Control, Academy of
Mathematics and Systems Science, Chinese Academy of Sciences,
Beijing 100190, China; Yangtze Center of Mathematics, Sichuan
University, Chengdu 610064, China; and BCAM-Basque Center for
Applied Mathematics, Bizkaia Technology Park, Building 500, E-48160,
Derio, Basque Country, Spain. {\small\it e-mail:} {\small\tt
xuzhang@amss.ac.cn}.}}
\maketitle
\begin{abstract}

This paper is addressed to the well-posedness of some linear and
semilinear backward stochastic differential equations with general
filtration, without using the Martingale Representation Theorem. The
point of our approach is to introduce a new notion of solution,
i.e., the transposition solution, which coincides with the usual
strong solution when the filtration is natural but it is more
flexible for the general filtration than the existing notion of
solutions. A comparison theorem for transposition solutions is also
presented.

\end{abstract}

\bs

\no{\bf 2010 Mathematics Subject Classification}.  Primary 60H10;
Secondary 34F05, 93E20.

\bs

\no{\bf Key Words} Backward stochastic differential equations,
transposition solution, filtration, comparison theorem.

\newpage

\section{Introduction}\label{s1}

Let $T>0$ and   $(\O,\cF,\dbF,\dbP)$ be a complete filtered
probability space with $\dbF=\{\cF_t\}_{t\in[0,T]}$, on which a
$1$-dimensional standard Brownian motion $\{w(t)\}_{t\in[0,T]}$ is
defined. We denote by $L_{\cF_t}^2(\O;\dbR^n)$  ($n\in \dbN$) the
Hilbert space consisting of all $\cF_t$-measurable ($\dbR^n$-valued)
random variables $\xi:\O\to \dbR^n$ such that
$\mathbb{E}|\xi|_{\dbR^n}^2 < \infty$, with the canonical inner
product; by $L^{2}_{\dbF}(\O;L^r(0,T;\dbR^n))$ ($1\leq r\leq
\infty$) the Banach space consisting of all $\dbR^n$-valued
$\{\cF_t\}$-adapted stochastic processes $X(\cdot)$ such that
$\mathbb{E}(|X(\cdot)|^2_{L^r(0,T;\dbR^n)}) < \infty$, with the
canonical norm; by $L^{2}_{\dbF}(\O;C([0,T];\dbR^n))$ the Banach
space consisting of all $\dbR^n$-valued $\{ {\cal F}_t \}$-adapted
continuous processes  $X(\cdot)$ such that
$\mathbb{E}(|X(\cdot)|^2_{L^{\infty}_{\dbF}(0,T;\dbR^n)}) < \infty$,
with the canonical norm; by $\cM_\dbF^2([0,T];\dbR^n)$ the Hilbert
space consisting of all $\dbR^n$-valued square integrable
$\{\cF_t\}$-martingales, with the canonical inner product; and by
$\cM_{0,\dbF}^2([0,T];\dbR^n)$ the closed subspace $\{X(\cd)\in
\cM_\dbF^2([0,T];\dbR^n)\;|\;X(0)=0\ \as\2n\}$ of
$\cM_\dbF^2([0,T];\dbR^n)$ with the inherited topology. Also, we
denote by $D([0,T];\dbR^n)$ the Banach space of all c\`adl\`ag
(i.e., right continuous with left limits) functions from $[0,T]$ to
$\dbR^n$, endowed with the inherited topology from
$L^{\infty}(0,T;\dbR^n)$ rather than the Skorokhod topology; and by
$L^{2}_{\dbF}(\O;D([0,T];\dbR^n))$ the Banach space consisting of
all $\dbR^n$-valued $\{ {\cal F}_t \}$-adapted c\`adl\`ag processes
$X(\cdot)$ such that
$\mathbb{E}(|X(\cdot)|^2_{L^{\infty}_{\dbF}(0,T;\dbR^n)}) < \infty$,
with the canonical norm. For any $t\in [0,T]$, one can define the
spaces $L^2_{\dbF}(\O;L^r(t,T;\dbR^n))$,
$L^2_{\dbF}(\O;C([t,T];\dbR^n))$, $L^2_{\dbF}(\O;D([t,T];\dbR^n))$
and so on in a similar way. Denote by $\lan \cd,\cd\ran$ the usual
scalar product in $\dbR^n$.

This paper is devoted to a study of the well-posedness for the
following semilinear backward stochastic differential equation (BSDE
for short)
\begin{eqnarray}\label{system1}
\left\{
\begin{array}{lll}
\ds dy (t) = f(t,y(t),Y(t))dt + Y(t) dw (t) & \mbox{ in } [0,T], \\
 \ns\ds y(T) = y_T,
\end{array}
\right.
 \end{eqnarray}
where $y_T\in L_{\cF_T}^2(\O;\dbR^n)$, $f(\cd,\cdot,\cdot)$
satisfies $f(\cd,0,0)\in L^2_{\dbF}(\O;L^1(0,T;\dbR^n))$ and, for
some constant $K>0$,
 \bel{Lm1}
 |f(t,p_1,q_1)-f(t,p_2,q_2)|\le
K(|p_1-p_2|+|q_1-q_2|),\q t\in [0,T], \hb{ a.s.},
\forall\;p_1,p_2,q_1,q_2\in \dbR^n.
 \ee
(Clearly, one can consider similarly the general case that the term
$Y(t) dw (t)$ in (\ref{system1}) is replaced by
$\left[g(t,y(t))+Y(t)\right] dw (t)$ provided that $g(\cd,0)\in
L^2_{\dbF}(\O;L^2(0,T;\dbR^n)) $ and $g(\cd,\cd)$ is globally
Lipschitz continuous with respect to its second argument).

The study of BSDEs is stimulated by the classical works
\cite{B,B1,PP}. Now, it is well-known that BSDEs and its various
variants play important and fundamental roles in Stochastic Control
(\cite{P,YZ}), Mathematical Finance (\cite{D,KPQ,T}), Probability
and Stochastic Analysis (\cite{P2}), Partial Differential Equations
(\cite{MY,P1,T}) and so on.

When $\dbF$ is equal to the natural filtration $\dbW$ (generated by
the Brownian motion $\{w(\cd)\}$ and augmented by all the
$\dbP$-null sets), the well-posedness of equation (\ref{system1}) is
well understood (\cite{PP}). In this case, by definition,
$(y(\cdot), Y(\cdot)) \in L^{2}_{\dbF}(\O;C([0,T];\dbR^n)) \t
L^2_{\dbF}(\O;L^2(0,T;\dbR^n))$ is said to be a (strong) solution to
equation (\ref{system1}) if
 \bel{7e7}
 y (t) = y_T-\int_t^Tf(s,y(s),Y(s))ds -\int_t^T Y(s) dw
 (s),\qq\forall\;t\in [0,T].
 \ee
Clearly, the first step to establish the well-posedness of the
semilinear equation (\ref{system1}) is to study the same problem but
for the following linear BSDE with a non-homonomous term $f(\cd)\in
L^2_{\dbF}(\O;L^1(0,T;$ $\dbR^n))$:
\begin{eqnarray}\label{sm1}
\left\{
\begin{array}{lll}
\ds dy (t) = f(t)dt + Y(t) dw(t),  & t\in [0,T), \\
 \ns\ds y(T) = y_T.
\end{array}
\right.
 \end{eqnarray}
The main idea in \cite{PP} for solving equation (\ref{system1}) with
$\dbF=\dbW$ is as follows: First, for (\ref{sm1}), noting that the
following process
 \bel{7e4}
 M(t)=\dbE\Big(y_T-\int_0^Tf(s)ds\;\Big|\;\cF_t\Big)
 \ee
is a $\{\cF_t\}$-martingale, and using the Martingale Representation
Theorem (valid only for the case $\dbF=\dbW$), one can find a
$Y(\cd)\in L^2_{\dbW}(\O;L^2(0,T;\dbR^n)) $ such that
 \bel{7e5}
 M(t)=M(0)+\int_0^tY(s)dw(s).
 \ee
Putting
 \bel{7e6} y(t)=M(t)+\int_0^tf(s)ds,
 \ee
one then finds the unique strong solution $(y(\cd),Y(\cd))\in
L^{2}_{\dbW}(\O;C([0,T];\dbR^n)) \t L^2_{\dbW}(\O;L^2(0,T;\dbR^n))$
for the linear BSDE (\ref{sm1}). Based on this and using the Picard
iteration argument, the desired well-posedness for equation
(\ref{system1}) follows.

It is easy to see that the Martingale Representation Theorem plays a
crucial role for the above mentioned well-posedness result for
equation (\ref{system1}) with natural filtration. In the general
case when the filtration $\dbF$ is not equal to the natural one,
$\dbW$ might be a proper sub-class of $\dbF$, and therefore, the
Martingale Representation Theorem fails. As far as we know, there
exists only a very few works addressing the well-posedness of
equation (\ref{system1}) with the general filtration (\cite{KH,
LLQ}).

The main idea  to study the well-posedness of BSDEs in \cite{KH} is
as follows. Consider first equation (\ref{sm1}). Since the
filtration $\dbF$ is not equal to the natural one, the following
 \begin{equation}\label{thi0rd1}
\cM_{0,\dbM,\dbF}^2([0,T];\dbR^n)\=\left.\left.\left\{\int_0^\cdot
g(s)dw(s)\;\right|\;\right.g(\cd)\in
 L^2_{\dbF}(\O;L^2(0,T;\dbR^n))\right\}
 \end{equation}
is a proper subspace of $\cM_{0,\dbF}^2([0,T];\dbR^n)$. Then one has
the following (unique) orthogonal decomposition:
\begin{equation}\label{third1}
 M(\cdot)-M(0)=P(\cd)+Q(\cd),
\end{equation}
for some $P(\cd)\in\cM_{0,\dbM,\dbF}^2([0,T];\dbR^n)$ and $Q(\cd)\in
\left(\cM_{0,\dbM,\dbF}^2([0,T];\dbR^n)\right)^\perp$. By
(\ref{thi0rd1}), there is a $Y(\cdot)\in
L^2_{\dbF}(\O;L^2(0,T;\dbR^n))$ such that
 \begin{equation}\label{thi1rd1}
  P(t) = \int_0^t
  {Y}(s)dw(s).
   \end{equation}
Still, we define $y(\cd)$ as in (\ref{7e6}). It is easy to check
that
 $
(y(\cd),Q(\cd),{Y}(\cd))\in L^{2}_{\dbF}(\O;D([0,T];\dbR^n)) \t
\left(\cM_{0,\dbM,\dbF}^2([0,T];\dbR^n)\right)^\perp\times
L^2_{\dbF}(\O;L^2(0,T;\dbR^n))
$
 is the unique solution of the
following equation
\begin{equation}\label{third2}
 y (t) = y_T+Q(t)-Q(T)-\int_t^T f(s)ds -\int_t^T Y(s) dw
 (s),\qq\forall\;t\in [0,T].
 \end{equation}
This means that (\ref{third2}) is another reasonable ``modification"
of the linear BSDE (\ref{sm1}) (by adding another corrected term
$Q(\cd)$). Similar to the above, by utilizing the Picard iteration
argument, one can study the well-posedness of equation
(\ref{system1}) (by adding one more corrected term $dQ(t)$ in the
right hand side of the first equation in (\ref{system1})). Note that
the appearance of this extra term $Q(\cd)$ makes the rigorous
analysis on the properties of $y(\cd)$ and $Y(\cd)$ much more
complicated than the case of natural filtration.  For example, one
needs to use some deep results in martingale theory (e.g.,
\cite[Chapter VIII]{DM}) to establish the duality relationship (like
(\ref{sol1}) below) between this sort of modified BSDEs and the
usual (forward) stochastic differential equations although it is not
difficult to give the desired relationship formally. Meanwhile, one
knows very little about $\cM_{0,\dbM,\dbF}^2([0,T];\dbR^n)$ (which
is actually introduced to replace the use of Martingale
Representation Theorem), and therefore, it seems very difficult to
``compute" the above $ Y(\cdot)$ in (\ref{thi1rd1}).

In \cite{LLQ}, the authors developed another approach to address the
well-posedness of BSDEs. The main idea in \cite{LLQ} for solving
equation (\ref{sm1}) (with general filtration) is as follows.
Although formula (\ref{7e5}) does not make sense any more,
$M(\cd)\in \cM_\dbF^2([0,T];\dbR^n)$ and $y(\cd)\in
L^{2}_{\dbF}(\O;D([0,T];\dbR^n))$\footnote{In \cite{LLQ}, the
authors asserted that $y(\cd)\in L^{2}_{\dbF}(\O;C([0,T];\dbR^n))$
(in terms of our notation). But it seems to us that this should be a
misprint.} are still well-defined respectively by (\ref{7e4}) and
(\ref{7e6}), and verifies $M(0)=y(0),\as$ Then, it is easy to check
that the above $(y(\cd),M(\cd))$ is the unique solution of the
following equation
  \bel{7e8}
 y (t) = y_T-\int_t^Tf(s)ds +M(t)-M(T),\qq\forall\;t\in [0,T]
 \ee
in the solution space
 \bel{9e2}
 \Upsilon\=\Big\{(h(\cd),N(\cd))\in
L^{2}_{\dbF}(\O;D([0,T];\dbR^n))\times
\cM_\dbF^2([0,T];\dbR^n)\;\Big|\;N(0)=h(0) \ \as\2n\Big\}.
 \ee
This means that (\ref{7e8}) is a reasonable ``modification" of the
linear BSDE (\ref{sm1}). Starting from this and using the Picard
iteration argument once more, one can study the well-posedness of
equation (\ref{system1}) (with a suitable modification) (See
\cite{LLQ} for more details). This approach does not need to use the
Martingale Representation Theorem, either. However, the adjusting
term $Y(\cd)$ in (\ref{sm1}) (or more generally, in (\ref{system1}))
is then suppressed. Note that this term plays a crucial role in some
problems, say the Pontryagin-type maximum principle for general
stochastic optimal control problems (\cite{P,YZ} and the references
therein). On the other hand, it seems to be very difficult to give
the duality analysis on solutions of equation (\ref{7e8}) (or the
modified version of  (\ref{system1})).

In this paper, we shall present a different approach to treat the
well-posedness of BSDEs with general filtration. Our idea is as
follows. Fixing $t\in [0,T]$, we consider the following linear
(forward) stochastic differential equation
\begin{eqnarray}\label{system2}
\left\{
\begin{array}{lll}
\ds dz(\tau)  = u(\tau)d\tau + v(\tau) dw (\tau), & \tau\in (t,T], \\
 \ns\ds z(t) = \eta.
\end{array}
\right.
 \end{eqnarray}
It is clear that, for given $u(\cdot)\in
L^2_{\dbF}(\O;L^1(t,T;\dbR^n))$, $v(\cdot)\in
L^2_{\dbF}(\O;L^2(t,T;\dbR^n))$ and $\eta\in
 L^2_{\cF_t}(\O;\dbR^n)$, equation (\ref{system2}) admits a unique strong solution $z(\cdot)\in
L^{2}_{\dbF}(\O;C([t,T];\dbR^n))$. Now, if equation (\ref{system1})
admits a strong solution $(y(\cdot), Y(\cdot)) \in
L^{2}_{\dbF}(\O;C([0,T];\dbR^n)) \t
L^{2}_{\dbF}(0,T;L^2(\O;\dbR^n))$ (say, when $\dbF=\dbW$), then,
applying It\^o's formula to $\lan z(t),y(t)\ran$, it is easy to
check that
 \bel{sol1}
 \ba{ll}
 \ds\mathbb{E}\lan z(T),y_T\ran - \mathbb{E}\lan
\eta,y(t)\ran\\
 \ns
\ds = \mathbb{E}\int_t^T \lan z(\tau),f(\tau,y(\tau),Y(\tau))\ran
d\tau + \mathbb{E}\int_t^T \lan u(\tau),y(\tau)\ran d\tau +
\mathbb{E} \int_t^T\lan v(\tau), Y(\tau)\ran d\tau.
 \ea
\ee
This inspires us to introduce the following new notion for the
solution of equation (\ref{system1}).

\begin{definition}\label{def of solution}
We call $(y(\cdot), Y(\cdot)) \in L^{2}_{\dbF}(\O;D([0,T];\dbR^n))
\t L^2_{\dbF}(\O;L^2(0,T;\dbR^n))$ a transposition solution of
equation (\ref{system1}) if for any $t\in [0,T]$, $u(\cdot)\in
L^2_{\dbF}(\O;L^1(t,T;\dbR^n))$, $v(\cdot)\in
 L^2_{\dbF}(\O;L^2(t,T;\dbR^n))$ and $\eta\in
 L^2_{\cF_t}(\O;\dbR^n)$,
identity (\ref{sol1}) holds.
\end{definition}

The main purpose of this paper is to show that equation
(\ref{system1}) is well-posed in the above transposition sense.
Clearly, any transposition solution of equation (\ref{system1})
coincides with its strong solution whenever the filtration $\dbF$ is
natural. Note that, in the general case, the space for the first
component of the solution is chosen to be
$L^{2}_{\dbF}(\O;D([0,T];\dbR^n))$ rather than
$L^{2}_{\dbF}(\O;C([0,T];\dbR^n))$. This is quite natural because
the filtration $\dbF$ is assumed only to be right-continuous.

Our approach is motivated by the classical transposition method in
solving the non-homogeneous boundary value problems for partial
differential equations  (\cite{LM}) and especially the boundary
controllability problem for hyperbolic equations (\cite{MR931277}).
On the other hand, one can find a rudiment of our approach at
\cite[pp. 353--354]{YZ} though the space for $y(\cdot)$ was chosen
to be $L^{2}_{\dbF}(\O;L^2(0,T;\dbR^n))$ and the filtration was
assumed to be natural there. The main advantage of our approach
consists in the fact that the duality analysis is contained in the
definition of solutions, and therefore, we do not need to utilize
the deep result in martingale theory to deduce this sort of duality
relationship any more, and one can easily deduce a similar
comparison theorem for transposition solutions of (\ref{system1}) by
using almost the same approach as in the case of natural filtration
(\cite{KPQ}). Also, it is even easier (and therefore we omit the
details) to establish a Pontryagin-type maximum principle for
general stochastic optimal control problems than to solve the same
problem with the natural filtration (\cite{P,YZ}) because, again,
the desired duality analysis is contained in the definition of
transposition solution. Moreover, by our method, the adjusting term
$Y(\cd)$ is obtained by the standard Riesz Representation Theorem
for Hilbert Space, and therefore, one can utilize the theory from
Hilbert Spaces to characterize $Y(\cdot)$, or even give a numerical
approach for $Y$ (see Remark \ref{numer}) although the detailed
analysis is beyond the scope of this paper.

People may be unsatisfied with our definition on the transposition
solution of (\ref{system1}) because one does not see what equation
this solution satisfies. However, starting from our transposition
solution of (\ref{system1}), one can obtain a corrected form of this
equation, i.e., equation (\ref{9e21}) in Section 4. Then, by
introducing suitably a corrected solution of (\ref{system1}) (See
Definition \ref{dn}), we obtain also a corresponding well-posedness
result (See Corollary \ref{th1}).

The rest of this paper is organized as follows. In Section 2, we
show some useful preliminary results. Section 3 is addressed to the
well-posedness of the linear BSDE (\ref{sm1}). Then, we prove the
well-posedness of the semilinear BSDE (\ref{system1}) in Section 4.
Finally, in Section 5, we present a comparison theorem for
transposition solutions of (\ref{system1}) in one dimension.

\section{Preliminaries}

In this section, we collect some preliminary results which will be
useful in the sequel.

Fix any $t_1$ and $t_2$ satisfying $0\leq t_2 \leq t_1 \leq T$.
First of all, we need the following Riesz-type Representation
Theorem, which is a special case of the known result in
\cite[Corollary 2.3 and  Remark 2.4]{LYZ}.

\bl\label{2l1}
For any $r\in [1,\infty)$, it holds that
 $$\left(L^2_\dbF(\O;L^r(t_2,t_1;\dbR^n))\right)^*=L^2_\dbF(\O;L^{r'}(t_2,t_1;\dbR^n)),
 $$
where $r'=r/(r-1)$ if $r\not=1$; $r'=\infty$ if $r=1$.
\el

Next, we need the following simple result (whose proof is direct,
and therefore we omit the details).

\bl\label{2l2}
There is a constant $C$, depending only on $T$, such that for any
$\big(u(\cd), v(\cd),\eta\big)\in L^2_{\dbF}(\O;L^1(t,T;\dbR^n))\t
L^2_{\dbF}(\O;L^2(t,T;\dbR^n))\times L^2_{\cF_t}(\O;\dbR^n)$, the
solution $z(\cdot)\in L^{2}_{\dbF}(\O;C([t,T];\dbR^n))$ of equation
(\ref{system2}) satisfies
  \bel{2.2}
  \ba{ll}\ds
 |z(\cdot)|_{L^{2}_{\dbF}(\O;C([t,T];\dbR^n))}\\
 \ns
 \ds\le
C\left|\big(u(\cd), v(\cd),\eta\big)\right|_{
L^2_{\dbF}(\O;L^2(t,T;\dbR^n))\t
L^2_{\dbF}(\O;L^2(t,T;\dbR^n))\times
L^2_{\cF_t}(\O;\dbR^n)},\qq\forall\;t\in [0,T].
 \ea
 \ee
\el

Further, we need the following result, which can be seen as a
variant of the classical Lebesgue Theorem (on Lebesgue point).

\bl\label{2l3}
Assume that $p\in(1,\infty]$, $q=\left\{\ba{ll}\frac{p}{p-1}&\hb{if
}\ p\in (1,\infty),\\[2mm] 1&\hb{if }\ p=\infty,\ea\right.$ $f_1\in
L^p_{\dbF}(0,T;L^2(\O;\dbR^n))$ and $f_2\in
L^q_{\dbF}(0,T;L^2(\O;\dbR^n))$. Then
  \bel{2.21}
 \lim_{h\to0}\frac{1}{h}\int_t^{t+h}\mathbb{E}\lan f_1(t),f_2(\tau)\ran d\tau=\mathbb{E}\lan f_1(t),f_2(t)\ran,\qq t\in [0,T]\ \ae
 \ee
\el

{\bf Proof.} We consider the case that $h\to0+$ (The case that
$h\to0-$ can be considered similarly). Let
$$\tilde{f_2} = \left\{\begin{array}{ll} \ds f_2,&  t\in [0,T]\\
\ns\ds 0, & t\in (T,2\,T]. \end{array}\right.$$ Obviously,
$\tilde{f_2}\in L^q_{\dbF}(0,2\,T;L^2(\O;\dbR^n))$ and
$$|\tilde{f_2}|_{L^q_{\dbF}(0,2\,T;L^2(\O;\dbR^n))}=|\tilde{f_2}|_{L^q_{\dbF}(0,T;L^2(\O;\dbR^n))}=|f_2|_{L^q_{\dbF}(0,T;L^2(\O;\dbR^n))}.$$ Since
$C([0,2\,T];L^2(\O;\dbR^n))$ is dense in
$L^q_{\dbF}(0,2T;L^2(\O;\dbR^n))$, for any $\e>0$, one can find
$f_2^0 \in C([0,2\,T];L^2(\O;\dbR^n))$ such that
 \bel{2-e1}
 |\tilde{f_2}-f_2^0|_{L^q_{\dbF}(0,2T;L^2(\O;\dbR^n))}\le
 \e.
 \ee
By the uniform continuity of $f_2^0(\cd)$ in $L^2(\O;\dbR^n)$, one
can find a $\d=\d(\e)>0$ such that
 \bel{2-e2}
 |f_2^0(s_1)-f_2^0(s_2)|_{L^2(\O;\dbR^n)}\le
 \e,\q \forall\;s_1,s_2\in [0,2\,T]\hb{ satisfying }|s_1-s_2|\le \d.
 \ee

Thanks to (\ref{2-e2}), we see that, when $h\le \d$, it holds that
 \bel{1e1}
 \ba{ll}
 \ds \q\int_0^T\left|\frac{1}{h}\int_t^{t+h}\mathbb{E}\lan f_1(t),f_2^0(\tau)\ran d\tau-\mathbb{E}\lan
 f_1(t),f_2^0(t)\ran\right|dt\\\ns
  \ds=\frac{1}{h}\int_0^T\left|\int_t^{t+h}\mathbb{E}\lan
 f_1(t),f_2^0(\tau)-f_2^0(t)\ran d\tau \right|dt\\\ns
  \ds\le \frac{1}{h}\int_0^T\int_t^{t+h}|f_1(t)|_{L^2(\O;\dbR^n)}|f_2^0(\tau)-f_2^0(t)|_{L^2(\O;\dbR^n)} d\tau dt\\\ns
  \ds\le \frac{\e}{h}\int_0^T\int_t^{t+h}|f_1(t)|_{L^2(\O;\dbR^n)} d\tau
  dt=\e\int_0^T|f_1(t)|_{L^2(\O;\dbR^n)} dt\le C\e |f_1|_{
L^p_{\dbF}(0,T;L^2(\O;\dbR^n))}.
 \ea
 \ee
Also, by (\ref{2-e1}), we have
 \bel{1e2}
 \ba{ll}
 \ds \q\int_0^T\left|\mathbb{E}\lan
 f_1(t),\tilde{f_2}(t)\ran-\mathbb{E}\lan
 f_1(t),f_2^0(t)\ran\right|dt\\\ns
  \ds\le |f_1|_{
L^p_{\dbF}(0,T;L^2(\O;\dbR^n))}|\tilde{f_2}-f_2^0|_{L^q_{\dbF}(0,2T;L^2(\O;\dbR^n))}\le
\e |f_1|_{ L^p_{\dbF}(0,T;L^2(\O;\dbR^n))}.
 \ea
 \ee
Further, using (\ref{2-e1}) again, we find
 \bel{1e3}
 \ba{ll}
 \ds \q\int_0^T\left|\frac{1}{h}\int_t^{t+h}\mathbb{E}\lan f_1(t),\tilde{f_2}(\tau)\ran d\tau-\frac{1}{h}\int_t^{t+h}\mathbb{E}\lan f_1(t),f_2^0(\tau)\ran d\tau\right|dt\\\ns
  \ds=\frac{1}{h}\int_0^T\left|\int_t^{t+h}\mathbb{E}\lan
 f_1(t),\tilde{f_2}(\tau)-f_2^0(\tau)\ran d\tau \right|dt\\\ns
  \ds\le \frac{1}{h}\int_0^T\int_t^{t+h}|f_1(t)|_{L^2(\O;\dbR^n)}|\tilde{f_2}(\tau)-f_2^0(\tau)|_{L^2(\O;\dbR^n)} d\tau dt\\\ns
  \ds\le \frac{1}{h}\left[\int_0^T\int_t^{t+h}|f_1(t)|_{L^2(\O;\dbR^n)}^pd\tau dt\right]^{1/p}\left[\int_0^T\int_t^{t+h}|\tilde{f_2}(\tau)-f_2^0(\tau)|
  _{L^2(\O;\dbR^n)}^q d\tau dt\right]^{1/q}\\\ns
  \ds=|f_1|_{
L^p_{\dbF}(0,T;L^2(\O;\dbR^n))}\left[
\frac{1}{h}\int_0^T\int_0^h|\tilde{f_2}(t+\tau)-f_2^0(t+\tau)|
  _{L^2(\O;\dbR^n)}^q d\tau dt\right]^{1/q}\\\ns
  \ds=|f_1|_{
L^p_{\dbF}(0,T;L^2(\O;\dbR^n))}\left[
\frac{1}{h}\int_0^h\int_\tau^{T+\tau}|\tilde{f_2}(t)-f_2^0(t)|
  _{L^2(\O;\dbR^n)}^q  dtd\tau\right]^{1/q}\\\ns
  \ds\le |f_1|_{
L^p_{\dbF}(0,T;L^2(\O;\dbR^n))}\left[
\frac{1}{h}\int_0^h\int_0^T|\tilde{f_2}(t)-f_2^0(t)|
  _{L^2(\O;\dbR^n)}^q  dtd\tau\right]^{1/q}\le \e |f_1|_{
L^p_{\dbF}(0,T;L^2(\O;\dbR^n))}.
 \ea
  \ee

Combining (\ref{1e1}), (\ref{1e2}) and (\ref{1e3}), we conclude that
 $$
\int_0^T\left|\frac{1}{h}\int_t^{t+h}\mathbb{E}\lan
f_1(t),\tilde{f_2}(\tau)\ran d\tau-\mathbb{E}\lan
f_1(t),\tilde{f_2}(t)\ran\right|dt\le C\e |f_1|_{
L^p_{\dbF}(0,T;L^2(\O;\dbR^n))}.
  $$
Therefore,
 $$
\lim_{h\to0}\int_0^T\left|\frac{1}{h}\int_t^{t+h}\mathbb{E}\lan
f_1(t),\tilde{f_2}(\tau)\ran d\tau-\mathbb{E}\lan
f_1(t),\tilde{f_2}(t)\ran\right|dt=0,
  $$
which means that  $$
 \lim_{h\to0}\frac{1}{h}\int_t^{t+h}\mathbb{E}\lan f_1(t),\tilde{f_2}(\tau)\ran d\tau=\mathbb{E}\lan f_1(t),\tilde{f_2}(t)\ran,\qq t\in [0,T]\ \ae
$$
By this and the definition of $\tilde{f_2}(\cdot)$, we conclude that
\begin{eqnarray*}
&\,&\lim_{h\to0}\frac{1}{h}\int_t^{t+h}\mathbb{E}\lan
f_1(t),f_2(\tau)\ran d\tau =
\lim_{h\to0}\frac{1}{h}\int_t^{t+h}\mathbb{E}\lan
f_1(t),\tilde{f_2}(\tau)\ran d\tau =\mathbb{E}\lan
f_1(t),\tilde{f_2}(t)\ran
\\\ns&\,& \qq\qq\qq\qq\qq\qq\q\;\;=\mathbb{E}\lan f_1(t),f_2(t)\ran,\qq t\in [0,T]\ \ae
\end{eqnarray*}
 This
completes the proof of Lemma \ref{2l3}.\endpf

\section{Well-posedness of linear non-homonomous BSDEs}

In this section, as a key step to study the well-posedness of the
semilinear BSDE (\ref{system1}), we consider first the same problem
but for equation (\ref{sm1}). We have the following result.

\begin{theorem}\label{the1}
For any $f(\cd)\in L^2_{\dbF}(\O;L^1(0,T;\dbR^n))$ and any $y_T \in
L^2_{\cF_T}(\O;\dbR^n)$, system (\ref{sm1}) admits a unique
transposition solution $(y(\cdot), Y(\cdot)) \in
L^{2}_{\dbF}(\O;D([0,T];\dbR^n)) \t L^{2}_{\dbF}(\O;L^2(0,T;\dbR^n))
 $
(in the sense of Definition \ref{def of solution}). Furthermore,
there is a constant $C$, depending only on $T$, such that
 \bel{2.1}
 \ba{ll}\ds
 |(y(\cdot), Y(\cdot))|_{
L^{2}_{\dbF}(\O;D([t,T];\dbR^n)) \t L^{2}_{\dbF}(\O;L^2(t,T;\dbR^n))}\\
\ns \ds\le C\left[
 |f(\cd)|_{ L^2_{\dbF}(\O;L^1(t,T;\dbR^n))}+|y_T|_{
L^2_{\cF_T}(\O;\dbR^n)}\right], \qq\forall\;t\in [0,T].
 \ea
 \ee
\end{theorem}

{\bf Proof.} We divide the proof into several steps.

\ms

{\bf Step 1.} We define a linear functional $\ell$ on
$L^2_{\dbF}(\O;L^1(t,T;\dbR^n))\t
L^2_{\dbF}(\O;L^2(t,T;\dbR^n))\times L^2_{\cF_t}(\O;\dbR^n)$ as
follows:
$$
 \ba{ll}
 \ds
\ell\big(u(\cd), v(\cd),\eta\big) = \mathbb{E}\lan z(T),y_T\ran -
\mathbb{E}\int_t^T \lan z(\tau),f(\tau)\ran dt,\\
 \ns
 \ds\qq\qq\qq\forall\; \big(u(\cd), v(\cd),\eta\big)\in
L^2_{\dbF}(\O;L^1(t,T;\dbR^n))\t
L^2_{\dbF}(\O;L^2(t,T;\dbR^n))\times L^2_{\cF_t}(\O;\dbR^n),
 \ea
$$
where $z(\cdot)\in L^{2}_{\dbF}(\O;C([t,T];\dbR^n))$ solves equation
(\ref{system2}).

Using the H\"older inequality and Lemma \ref{2l2},  it is easy to
show that
 \bel{l bound}
 \ba{ll}
 \ds \left|\ell\big(u(\cd),
v(\cd),\eta\big)\right| \\
 \ns
\ds\leq |z(T)|_{ L^2_{\cF_T}(\O;\dbR^n)}|y_T|_{
L^2_{\cF_T}(\O;\dbR^n)} +
|z(\cdot)|_{L^{2}_{\dbF}(\O;C([t,T];\dbR^n))}
|f|_{L^2_{\dbF}(\O;L^2(t,T;\dbR^n))}\\
 \ns
\ds\leq C\left[
 |f(\cd)|_{ L^2_{\dbF}(\O;L^1(t,T;\dbR^n))}+|y_T|_{
L^2_{\cF_T}(\O;\dbR^n)}\right]\\
 \ns
\ds\qq\times \left|\big(u(\cd), v(\cd),\eta\big)\right|_{
L^2_{\dbF}(\O;L^1(t,T;\dbR^n))\t
L^2_{\dbF}(\O;L^2(t,T;\dbR^n))\times L^2_{\cF_t}(\O;\dbR^n)},
\q\forall\;t\in [0,T],
 \ea\ee
where $C=C(T)$ is independent of $t$. From (\ref{l bound}), we know
$\ell$ is a bounded linear functional on
$L^2_{\dbF}(\O;L^2(t,T;\dbR^n))\t
L^2_{\dbF}(\O;L^2(t,T;\dbR^n))\times L^2_{\cF_t}(\O;\dbR^n)$. Now,
by means of Lemma \ref{2l1}, we conclude that there exist
$y^t(\cdot)\in L^2_{\dbF}(\O;L^\infty(t,T;\dbR^n))$, $Y^t(\cd) \in
L^2_{\dbF}(\O;L^2(t,T;\dbR^n))$ and $\varsigma^t \in
L^2_{\cF_t}(\O;\dbR^n)$ such that
 \bel{eq1}
 \ba{ll}
 \ds\mathbb{E}\lan z(T),y_T\ran -
\mathbb{E}\int_t^T \lan z(\tau),f(\tau)\ran d\tau\\
 \ns\ds
 =  \mathbb{E}\int_t^T
\lan u(\tau),y^t(\tau)\ran d\tau + \mathbb{E} \int_t^T\lan
v(\tau),Y^t(\tau)\ran d\tau +\mathbb{E} \lan \eta,\varsigma^t\ran.
 \ea
\ee
It is clear that $\varsigma^T=y_T$. Furthermore, there is a positive
constant $C=C(T)$, independent of $t$, such that
 \bel{2.3}
 \ba{ll}\ds
|(y^t(\cdot), Y^t(\cd),\varsigma^t)|_{
L^2_{\dbF}(\O;L^\infty(t,T;\dbR^n)) \t L^2_{\dbF}(\O;L^2(t,T;\dbR^n))\times L^2_{\cF_t}(\O;\dbR^n)}\\
\ns \ds\le C\left[
 |f(\cd)|_{ L^2_{\dbF}(\O;L^1(t,T;\dbR^n))}+|y_T|_{
L^2_{\cF_T}(\O;\dbR^n)}\right], \qq\forall\;t\in [0,T].
 \ea
 \ee

{\bf Step 2.} Note that the ``solution" $(y^t(\cdot), Y^t(\cd))$
(obtained in Step 1) may depend on $t$. In  this step, we shall show
the time consistency of  $(y^t(\cdot), Y^t(\cd))$, i.e., for any
$t_1$ and $t_2$ satisfying $0\leq t_2 \leq t_1 \leq T$, it holds
that
 \bel{2e4}
 \big(y^{t_2} (\tau,\o),Y^{t_2} (\tau,\o)\big)=\big( y^{t_1}(\tau,\o), Y^{t_1}(\tau,\o)\big),\qq (\tau,\o) \in [t_1,T]\t\O\ \ \ae
 \ee

Note that the solution $z(\cd)$ of equation (\ref{system2}) depends
on $t$, and therefore, we also denote it by $z^t(\cd)$ (whenever
there exists a possible confusion). To show $y^{t_2} (\tau,\o)=
y^{t_1}(\tau,\o)$ for a.e. $(\tau,\o) \in [t_1,T]\t\O$, we fix any
$\varrho (\cd)\in L^2_{\dbF}(\O;L^1(t_1,T;\dbR^n))$ and choose
$t=t_1$, $\eta = 0$, $v(\cdot) = 0$ and $u(\cdot)=\varrho (\cd)$ in
equation (\ref{system2}). From (\ref{eq1}), we see that
\begin{equation}\label{eq2}
 \mathbb{E}\lan z^{t_1}(T),y_T\ran -
\mathbb{E}\int_{t_1}^T \lan z^{t_1}(\tau),f(\tau)\ran d\tau
 =  \mathbb{E}\int_{t_1}^T
\lan \varrho (\tau),y^{t_1}(\tau)\ran d\tau.
\end{equation}
On the other hand, choosing $t=t_2$, $\eta = 0$, $v(\cdot) = 0$ and
$u(t,\omega) = \chi_{[t_1,T]}(t) \varrho (t,\omega)$ in equation
(\ref{system2}). It is clear that $$z^{t_2}(\cdot) =\left\{
\begin{array}{ll}\ds  z^{t_1}(\cdot),& t \in [t_1,T],\\
\ns\ds 0, & t \in [t_2,t_1). \end{array}\right.
 $$
In this case, by (\ref{eq1}), we have
\begin{equation}\label{eq3}
 \mathbb{E}\lan z^{t_1}(T),y_T\ran -
\mathbb{E}\int_{t_1}^T \lan z^{t_1}(\tau),f(\tau)\ran d\tau
 =  \mathbb{E}\int_{t_1}^T
\lan \varrho (\tau),y^{t_2}(\tau)\ran d\tau.
\end{equation}
From (\ref{eq2}) and (\ref{eq3}), we conclude that
 $$
 \mathbb{E}\int_{t_1}^T
\lan \varrho (\tau),y^{t_1}(\tau)\ran d\tau = \mathbb{E}\int_{t_1}^T
\lan \varrho (\tau),y^{t_2}(\tau)\ran d\tau,\qq \forall\; \varrho
(\cd)\in L^2_{\dbF}(\O;L^1(t_1,T;\dbR^n)).
 $$
From this, we see that $y^{t_2} (\tau,\o)= y^{t_1}(\tau,\o)$ for
$(\tau,\o) \in [t_1,T]\t\O\ \ \ae$

To show $Y^{t_2} (\tau,\o)= Y^{t_1}(\tau,\o)$ for a.e. $(\tau,\o)
\in [t_1,T]\t\O$, we fix any $\varsigma (\cd)\in
L^2_{\dbF}(\O;L^2(t_1,T;\dbR^n))$ and choose $t=t_1$, $\eta = 0$,
$u(\cdot) = 0$ and $v(\cdot)=\varsigma (\cd)$ in equation
(\ref{system2}) (and denote by $\bar z^{t_1}(\cd)$ the corresponding
solution of (\ref{system2})). From (\ref{eq1}), we see that
\begin{equation}\label{e0q2}
 \mathbb{E}\lan \bar z^{t_1}(T),y_T\ran -
\mathbb{E}\int_{t_1}^T \lan \bar z^{t_1}(\tau),f(\tau)\ran d\tau
 =  \mathbb{E}\int_{t_1}^T
\lan \varsigma (\tau),Y^{t_1}(\tau)\ran d\tau.
\end{equation}
On the other hand, choosing $t=t_2$, $\eta = 0$, $u(\cdot) = 0$ and
$v(t,\omega) = \chi_{[t_1,T]}(t) \varsigma (t,\omega)$ in equation
(\ref{system2})  (and denote by $\bar z^{t_2}(\cd)$ the
corresponding solution of (\ref{system2})). It is clear that
 $$\bar z^{t_2}(\cdot) =\left\{
\begin{array}{ll}\ds \bar  z^{t_1}(\cdot),& t \in [t_1,T],\\
\ns\ds 0, & t \in [t_2,t_1). \end{array}\right.
 $$
In this case, by (\ref{eq1}), we have
\begin{equation}\label{e-q3}
 \mathbb{E}\lan \bar z^{t_1}(T),y_T\ran -
\mathbb{E}\int_{t_1}^T \lan\bar  z^{t_1}(\tau),f(\tau)\ran d\tau
 =  \mathbb{E}\int_{t_1}^T
\lan \varsigma (\tau),Y^{t_2}(\tau)\ran d\tau.
\end{equation}
From (\ref{e0q2}) and (\ref{e-q3}), we conclude that
 $$
 \mathbb{E}\int_{t_1}^T
\lan \varsigma (\tau),Y^{t_1}(\tau)\ran d\tau =
\mathbb{E}\int_{t_1}^T \lan \varsigma (\tau),Y^{t_2}(\tau)\ran
d\tau,\qq \forall\; \varsigma (\cd)\in
L^2_{\dbF}(\O;L^2(t_1,T;\dbR^n)).
 $$
From this, we see that $Y^{t_2} (\tau,\o)= Y^{t_1}(\tau,\o)$ for
$(\tau,\o) \in [t_1,T]\t\O\ \ \ae$ Hence, (\ref{2e4}) is verified.

Put
 \bel{2e10}
 y(t,\o)=y^0(t,\o),\qq Y (t,\o)= Y^0(t,\o),\qq \forall\;(t,\o) \in
[0,T]\t\O.
 \ee
Then, in view of (\ref{2e4}), it follows that
 \bel{2e11}
 \big(y^t (\tau,\o),Y^t
(\tau,\o)\big)=\big( y(\tau,\o),Y (\tau,\o)\big), \qq (\tau,\o) \in
[t,T]\t\O\ \ \ae
 \ee
Combining (\ref{eq1}) and (\ref{2e11}), we find that
 \bel{eq16}
 \ba{ll}
 \ds\q\mathbb{E}\lan z(T),y_T\ran - \mathbb{E} \lan \eta,\varsigma^t\ran\\
 \ns\ds
 =\mathbb{E}\int_t^T \lan
z(\tau),f(\tau)\ran d\tau+ \mathbb{E}\int_t^T \lan
u(\tau),y(\tau)\ran d\tau +\mathbb{E} \int_t^T\lan
v(\tau),Y(\tau)\ran d\tau .
 \ea
\ee

\ms

{\bf Step 3.} We show in this step that $\varsigma^t$ has a
c\`adl\`ag modification. For this, clearly, it suffices to show that
 \bel{X}
 X(t)\=\varsigma^t - \int_0^t f(s)ds,\qq t\in[0,T]
 \ee
is a $\{\cF_t\}$-martingale. The rest of this step is to show that
$\{ X(t)\}$ is a $\{\cF_t\}$-martingale.

First of all, we claim that, for each $t\in [0,T]$,
 \bel{6e1} \dbE\Big(y_T -
\int_t^T f(s)ds \;\Big|\;\cF_t\Big) = \varsigma^t,\ \ \as
 \ee
To show this, choosing $z(t)=\varsigma^t$, $u=0$ and $v=0$ in
(\ref{system2}), it follows that
$$
\dbE\lan\varsigma^t,y_T\ran - \dbE |\varsigma^t|^2 = \dbE\int_t^T
\lan \varsigma^t,f(s)\ran ds.
$$
This gives
\begin{equation}\label{M1}
\dbE \lan\varsigma^t,\dbE(y_T\;|\;\cF_t)\ran- \dbE |\varsigma^t|^2 =
\dbE  \Big\langle \varsigma^t,\dbE\Big(\int_t^T
f(s)ds\;\Big|\;\cF_t\Big) \Big\rangle.
\end{equation}
From equality (\ref{M1}), we have
\begin{equation}\label{M2}
\dbE \Big\langle \varsigma^t,\dbE\Big(y_T-\int_t^T
f(s)ds\;\Big|\;\cF_t\Big)\Big\rangle= \dbE |\varsigma^t|^2 .
\end{equation}
On the other hand, choosing  $z(t)=\dbE\Big(y_T-\int_t^T
f(s)ds\;\Big|\;\cF_t\Big)$, $u=0$ and $v=0$ in (\ref{system2}), we
obtain that
\begin{eqnarray}\label{M3}
&\,&\dbE \Big\langle\dbE\Big(y_T - \int_t^T f(s)ds\;\Big|\;\cF_t
\Big),y_T\Big\rangle -
\dbE \Big\langle \varsigma^t,\dbE\Big(y_T - \int_t^T f(s)ds \;\Big|\;\cF_t\Big)\Big\rangle\nonumber\\
&\,& = \dbE \Big\langle  \dbE\Big(y_T - \int_t^T f(s)ds \;\Big|\;
\cF_t\Big), \int_t^T f(s)ds \Big\rangle\nonumber\\
&\,& = \dbE \Big\langle \dbE\Big(y_T - \int_t^T f(s)ds \;\Big|\;
\cF_t\Big),\dbE\Big(\int_t^T f(s)ds\;\Big|\;\cF_t\Big) \Big\rangle.
\end{eqnarray}
From equality (\ref{M3}), we arrive at
\begin{equation}\label{M4}
\dbE\left|\dbE \Big(y_T - \int_t^T f(s)ds
\;\Big|\;\cF_t\Big)\right|^2 - \dbE
\Big\langle\varsigma^t,\dbE\Big(y_T-\int_t^T
f(s)ds\;\Big|\;\cF_t\Big)\Big\rangle = 0.
\end{equation}
Combining equality (\ref{M2}) and (\ref{M4}), we  end up with
\begin{equation*}
\dbE\left|\dbE\Big(y_T - \int_t^T f(s)ds \;\Big|\;\cF_t\Big) -
\varsigma^t\right|^2 = 0,
\end{equation*}
which gives (\ref{6e1}).

Next, combining (\ref{2.3}) and (\ref{6e1}), it is easy to see that
$\varsigma^\cd\in L^2_{\dbF}(\O;L^\infty(0,T;\dbR^n))$. Hence, $X
(\cd)\in L^2_{\dbF}(\O;L^2(0,T;\dbR^n))$.

Now, for any $\tau_1, \tau_2 \in [0,T]$ with $\tau_1 \leq \tau_2$,
by (\ref{6e1}), it follows that
\begin{eqnarray}
&\,&\dbE(X(\tau_2)\;|\;\cF_{\tau_1}) =
\dbE\Big(\varsigma^{\tau_2}-\int_0^{\tau_2}
f(s)ds\;\Big|\;\cF_{\tau_1}\Big)\nonumber\\
&\,& \qq\qq \qq\ \ = \dbE\left.\left[\dbE\Big(y_T - \int_{\tau_2}^T
f(s)ds\;\Big|\;\cF_{\tau_2}\Big)-\int_0^{\tau_2}
f(s)ds\;\right|\;\cF_{\tau_1}\right]\nonumber\\
&\,& \qq\qq \qq\ \  = \dbE\Big(y_T-\int_0^{T}
f(s)ds\;\Big|\;\cF_{\tau_1}\Big)\nonumber\\
&\,& \qq\qq \qq\ \ = \dbE\Big(y_T-\int_{\tau_1}^{T}
f(s)ds\;\Big|\;\cF_{\tau_1}\Big)-\int_0^{\tau_1}f(s)ds\nonumber\\
&\,& \qq\qq \qq\ \ = \varsigma^{\tau_1}- \int_0^{\tau_1}f(s)ds \nonumber\\
&\,& \qq\qq \qq\ \ = X(\tau_1), \q \mbox{ a.s. }
\end{eqnarray}
Therefore, $\{X(t)\}_{0\leq t\leq T}$ is a $\cF_{t}$-martingale.

\ms

{\bf Step 4.} In this step, we show that, for a.e $t\in [0,T]$,
 \bel{6e3}
 \varsigma^t= y(t)\ \ \as
 \ee
Fix any $\gamma \in L^2_{\cF_{t_2}}(\O;\dbR^n)$. Choosing $t=t_2$,
$u(\cd) = 0$, $v(\cd) = 0$ and $\eta = (t_1-t_2)\gamma$ in
(\ref{system2}), using (\ref{eq16}), we obtain that
\begin{equation}\label{eq6}
\mathbb{E}\lan (t_1-t_2)\gamma, y_T \ran-\mathbb{E}\lan
(t_1-t_2)\gamma, \varsigma^{t_2}\ran= \mathbb{E}\int_{t_2}^T \lan
(t_1-t_2)\gamma, f(\tau)\ran d\tau .
\end{equation}
Choosing $t=t_2$, $u(\tau,\o) = \chi_{[t_2,t_1]}(\tau)\gamma(\o)$,
$v(\cd) = 0$ and $\eta = 0$ in (\ref{system2}), using (\ref{eq16})
once more, we conclude that
\begin{equation}\label{eq7}
 \ba{ll}\ds
\mathbb{E}\lan (t_1-t_2)\gamma, y_T \ran\\
\ns \ds= \mathbb{E}\int_{t_2}^{t_1}\lan (\tau-t_2)\gamma,
f(\tau)\ran d\tau + \mathbb{E}\int_{t_1}^T \lan (t_1-t_2)\gamma,
f(\tau)\ran d\tau + \mathbb{E}\int_{t_2}^{t_1}\lan
\gamma,y(\tau)\ran d\tau.
 \ea
\end{equation}
From (\ref{eq6}) and (\ref{eq7}), we end up with
 $$
\mathbb{E}\lan \gamma, \varsigma^{t_2}\ran=
\frac{1}{t_1-t_2}\mathbb{E}\int_{t_2}^{t_1}\lan (\tau-t_2)\gamma,
f(\tau)\ran d\tau-\int_{t_2}^{t_1}\lan \gamma, f(\tau)\ran d\tau
+\frac{1}{t_1-t_2}\int_{t_2}^{t_1}\mathbb{E}\lan \gamma,y(\tau) \ran
d\tau.
 $$
It is easy to show that
 $$
 \lim_{t_1\to t_2+0} \frac{1}{t_1-t_2}\mathbb{E}\int_{t_2}^{t_1}\lan (\tau-t_2)\gamma,
f(\tau)\ran d\tau=\lim_{t_1\to t_2+0} \int_{t_2}^{t_1}\lan \gamma,
f(\tau)\ran d\tau=0,\qq\forall\;\gamma \in
L^2_{\cF_{t_2}}(\O;\dbR^n).
 $$
Hence,
\begin{equation}\label{eq8}
\lim_{t_1\to t_2+0} \frac{1}{t_1-t_2}\int_{t_2}^{t_1}\mathbb{E}\lan
\gamma,y(\tau) \ran d\tau=\mathbb{E}\lan \gamma,
\varsigma^{t_2}\ran,\qq\forall\;\gamma \in
L^2_{\cF_{t_2}}(\O;\dbR^n).
\end{equation}

Now, we need to compute the limit $\ds\lim_{t_1\to t_2+0}
\frac{1}{t_1-t_2}\int_{t_2}^{t_1}\mathbb{E}\lan \gamma,y(\tau) \ran
d\tau$ for some special $\gamma$. We consider first the simple case
that $L^2_{\cF_T}(\O;\dbR^n)$ is a separable Hilbert space. In this
case, one can find a sequence $\{\gamma_k\}_{k=1}^\infty$ which is
dense in $L^2_{\cF_T}(\O;\dbR^n)$. For each $k$, by the classical
Lebesgue Theorem (on Lebesgue point), we conclude that there is a
Lebesgue null set $E_k$ such that
 \bel{3e01}\ds\lim_{t_1\to t_2+0}
\frac{1}{t_1-t_2}\int_{t_2}^{t_1}\mathbb{E}\lan \gamma_k,y(\tau)
\ran d\tau=\lan \gamma_k,y(t_2)\ran, \qq\forall\; t_2\in
[0,T]\setminus E_k.
 \ee
Put $\ds E=\bigcup_{k=1}^\infty E_k$, whose Lebesque measure is $0$.
By (\ref{3e01}) and noting the density of
$\{\gamma_k\}_{k=1}^\infty$ in $L^2_{\cF_T}(\O;\dbR^n)$, it follows
that
 \bel{3e02}\ds\lim_{t_1\to t_2+0}
\frac{1}{t_1-t_2}\int_{t_2}^{t_1}\mathbb{E}\lan \gamma,y(\tau) \ran
d\tau=\lan \gamma,y(t_2)\ran, \qq\forall\;\gamma \in
L^2_{\cF_{t_2}}(\O;\dbR^n), \forall\; t_2\in [0,T]\setminus E.
 \ee
Combining (\ref{eq8}) and (\ref{3e02}), we find that $\mathbb{E}\lan
\gamma, \varsigma^{t_2}\ran=\mathbb{E}\lan \gamma,y(t_2)\ran$ for
any $\gamma \in L^2_{\cF_{t_2}}(\O;\dbR^n)$ and any $t_2\in
[0,T]\setminus E$. Hence, $\varsigma^t= y(t)$ in $[0,T]\t\O$, \ae

Now, we analyze the general case that $L^2_{\cF_T}(\O;\dbR^n)$ may
not be a separable Hilbert space. In this case, by (\ref{eq8}), we
conclude that
\begin{equation}\label{eq10}
\lim_{t_1\to t_2+0} \frac{1}{t_1-t_2}\int_{t_2}^{t_1}\mathbb{E}\lan
\varsigma^{t_2}-y(t_2),y(\tau) \ran d\tau=\mathbb{E}\lan
\varsigma^{t_2}-y(t_2), \varsigma^{t_2}\ran.
\end{equation}
Using Lemma \ref{2l3}, it follows
\begin{equation}\label{eq12}
\lim_{t_1\to t_2+0} \frac{1}{t_1-t_2}\int_{t_2}^{t_1}\mathbb{E}\lan
\varsigma^{t_2}-y(t_2),y(\tau) \ran d\tau=\mathbb{E}\lan
\varsigma^{t_2}-y(t_2),y(t_2)\ran,\qq t_2\in [0,T]\ \ae
\end{equation}
By (\ref{eq10})--(\ref{eq12}), we arrive at
\begin{equation}\label{eq14}
\mathbb{E}\lan \varsigma^{t_2}-y(t_2),
\varsigma^{t_2}\ran=\mathbb{E}\lan
\varsigma^{t_2}-y(t_2),y(t_2)\ran,\qq t_2\in [0,T]\ \ae
\end{equation}
By (\ref{eq14}), we find that $\mathbb{E}\left|
\varsigma^{t_2}-y(t_2)\right|^2=0$ for $ t_2\in [0,T]$ a.e., which
implies (\ref{6e3}) immediately.

Finally, combining (\ref{6e3}) and the result in Step 3 that
$\varsigma^t$ has a c\`adl\`ag modification, we see that there is a
c\`adl\`ag $\dbR^n$-valued process $\{\tilde y(t)\}_{t\in [0,T]}$
such that $y(\cd)=\tilde y$ in $[0,T]\times \Omega$ \ae It is easy
to check that $(\tilde y(\cdot),Y(\cdot))$ is a transposition
solution to equation (\ref{sm1}). To simplify the notation, we still
use $y$ instead of $\tilde y$ to denote the first component of the
solution. This means that equation (\ref{sm1}) admits one and only
one transposition solution $(y(\cdot), Y(\cdot)) \in
L^{2}_{\dbF}(\O;D([0,T];\dbR^n)) \t L^{2}_{\dbF}(\O;L^2(0,T;\dbR^n))
 $, which completes the proof of Theorem \ref{the1}.
\endpf

\br\label{3r2}
For the linear BSDE (\ref{sm1}), we may introduce another notion of
solution. We call $y(\cdot) \in L^{2}_{\dbF}(\O;D([0,T];\dbR^n))$ a
transposition pseudo-solution of equation (\ref{sm1}) if for any
$t\in [0,T]$, $u(\cdot)\in L^2_{\dbF}(\O;L^1(t,T;\dbR^n))$ and
$\eta\in
 L^2_{\cF_t}(\O;\dbR^n)$,
the following identity holds
 \bel{so-l1}
 \ba{ll}
 \ds\q\mathbb{E}\lan z(T),y_T\ran - \mathbb{E}\lan
\eta,y(t)\ran = \mathbb{E}\int_t^T \lan z(\tau),f(\tau)\ran d\tau +
\mathbb{E}\int_t^T \lan u(\tau),y(\tau)\ran d\tau,
 \ea
 \ee
where $z(\cdot)\in L^{2}_{\dbF}(\O;C([t,T];\dbR^n))$ solves equation
(\ref{system2}) with $v(\cdot)= 0$. Using almost the same proof as
that of Theorem \ref{the1}, one can show that, for any $f(\cd)\in
L^2_{\dbF}(\O;L^1(0,T;\dbR^n))$ and any $y_T \in
L^2_{\cF_T}(\O;\dbR^n)$, system (\ref{sm1}) admits a unique
transposition pseudo-solution $y(\cdot) \in
L^{2}_{\dbF}(\O;D([0,T];\dbR^n))$. Furthermore, there is a constant
$C$, depending only on $T$, such that
 $$
 |y(\cdot)|_{
L^{2}_{\dbF}(\O;D([t,T];\dbR^n)) }\le C\left[
 |f(\cd)|_{ L^2_{\dbF}(\O;L^1(t,T;\dbR^n))}+|y_T|_{
L^2_{\cF_T}(\O;\dbR^n)}\right], \qq\forall\;t\in [0,T].
 $$
It is clear that the transposition pseudo-solution $y(\cdot)$ of
equation (\ref{sm1}) coincides with the first component of the  the
transposition solution $y(\cdot)$ of equation (\ref{sm1}).
Nevertheless, the transposition pseudo-solution is not a good notion
for solution of equation (\ref{sm1}) because it does not reproduce
the strong solution even if the filtration is natural.
\er

\begin{remark}\label{numer}
At least conceptually, we can give a ``numerical" approach for BSDEs
with the general filtration in terms of the transposition solution.
Indeed, let $\{H_m\}_{m=1}^{+\infty}$ be a sequence of subspaces of
$L^{2}_{\dbF}(\O;L^2(0,T;\dbR^n))$ such that for any $g(\cdot)\in
L^{2}_{\dbF}(\O;L^2(0,T;\dbR^n))$, there exists a sequence
$\{g_m(\cd)\}_{m=1}^{+\infty}$ satisfies that
$$g_m(\cd)\in H_m \mbox{ and } \lim_{m\to{+\infty}}|g_m -
g|_{L^{2}_{\dbF}(\O;L^2(0,T;\dbR^n))}=0.
$$
Now, for a fixed $m\in\dbN$, choosing $t=0$, $\eta = 0$, $u(\cdot) =
0$ and $v(\cdot)=v_m(\cdot)\in H_m$ in equation (\ref{system2}) (and
denote by $\bar z_m(\cd)$ the corresponding solution of
(\ref{system2})).  From (\ref{eq1}), we see that
\begin{equation}\label{numer1}
\dbE\lan \bar z_m(T),y_T\ran - \mathbb{E}\int_0^T \lan \bar
z_m(\tau),f(\tau)\ran d\tau = \mathbb{E}\int_0^T \lan
v_m(\tau),Y(\tau)\ran d\tau.
\end{equation}
On the other hand, using the same argument to obtain $Y(\cdot)$ (by
Riesz's Representation Theorem), we can find a $Y_m(\cdot)\in H_m$
such that
\begin{equation}\label{numer2}
\dbE\lan \bar z_m(T),y_T\ran - \mathbb{E}\int_0^T \lan \bar
z_m(\tau),f(\tau)\ran d\tau = \mathbb{E}\int_0^T \lan
v_m(\tau),Y_m(\tau)\ran d\tau.
\end{equation}
This, together with \eqref{numer1}, implies that $Y_m(\cd)$ is the
orthogonal projection of $Y(\cd)$ to $H_m$. By the definition of
$H_m$, we know that $$
\lim_{m\to{+\infty}}|Y_m(\cdot)-Y(\cd)|_{L^{2}_{\dbF}(\O;L^2(0,T;\dbR^n))}=0.
$$ Therefore, one can get a ``good" approximation of $Y(\cd)$ if one
can choose a suitable sequence $\{H_m\}_{m=1}^{+\infty}$ such that
$Y_m(\cdot)$ (say, belongs to a finite dimensional space) can be
computed efficiently and that $Y_m(\cd)$ converges to $Y(\cd)$ in
$L^{2}_{\dbF}(\O;L^2(0,T;\dbR^n))$ in some sense. This will be done
in our forthcoming work.
\end{remark}

\br\label{3r1}
It is clear that one of the key observation in the proof of Theorem
\ref{the1} is that the process $ X(t)$ defined by (\ref{X}) is a
$\cF_t$-martingale. Combining this fact and (\ref{6e3}), we see that
the following process
 \bel{9e1}
 M(t)\=y(t)-\int_0^tf(s)ds,\qq\forall\; t\in [0,T]
  \ee
is a $\cF_t$-martingale as well. From this, it is easy to check that
$(y(\cd),M(\cd))$ is the unique solution of equation (\ref{7e8}) in
the solution space $\Upsilon$ (defined by (\ref{9e2})). Hence,
starting from our transposition solution $(y(\cd),Y(\cd))$ for the
linear BSDE (\ref{sm1}), one can re-construct the solution
$(y(\cd),M(\cd))$ introduced in \cite{LLQ}, through the relationship
(\ref{9e1}) between $M(\cd)$ and $y(\cd)$. Note however that one
cannot do the reverse because for the later one needs to represent
$Y(\cd)$ in terms of $M(\cd)$, which is exactly the concern of the
Martingale Representation Theorem. Nevertheless, since the solution
$(y(\cd),M(\cd))$ of equation (\ref{7e8}) is unique in $\Upsilon$,
it is easy to see that the first component of this solution
coincides with the first component of the transposition solution
$(y(\cd),Y(\cd))$ for equation (\ref{sm1}).
\er

\br
From the proof of Theorem \ref{the1}, it is easy to see why we
choose the space for the first component of the transposition
solution to be $L^{2}_{\dbF}(\O;D([0,T];\dbR^n))$ rather than $
L^{2}_{\dbF}(\O;C([0,T];\dbR^n))$ because a $\cF_t$-martingale has
only a  c\`adl\`ag modification. Indeed, as far as we know, there is
no general solution on the problem: Under what conditions, a
martingale has a continuous modification? We refer to \cite[Theorem
2.1.44]{YJ} for partial solution to this problem.
\er

Before ending this section, we put
 $$
 \ba{ll}
 L_{0,\cF_T}^2(\O;\dbR^n)\=\left.\left.\left\{h\in L_{\cF_T}^2(\O;\dbR^n)\;\right|\;\right.\dbE h=0\right\},\\
 \ns
\ds
L_{\dbW,\cF_T}^2(\O;\dbR^n)\=\left.\left.\left\{\int_0^Tg(s)dw(s)\;\right|\;\right.g(\cd)\in
 L^2_{\dbF}(\O;L^2(0,T;\dbR^n))\right\}.
 \ea
 $$
Clearly, $ L_{\dbW,\cF_T}^2(\O;\dbR^n)$ is a closed subspace of
$L_{0,\cF_T}^2(\O;\dbR^n)$. Generally, $
L_{\dbW,\cF_T}^2(\O;\dbR^n)$ is a proper subspace of $
L_{0,\cF_T}^2(\O;\dbR^n)$. Hence, the orthogonal complement space $
\left( L_{\dbW,\cF_T}^2(\O;\dbR^n)\right)^\perp$ is well-defined. We
have the following result.

\bp\label{3p1}
i) If $\ds y_T-\int_0^Tf(s)ds-\dbE\left(y_T-\int_0^Tf(s)ds\right)\in
L_{\dbW,\cF_T}^2(\O;\dbR^n)$, then the transposition solution
$(y(\cdot), Y(\cdot)) \in L^{2}_{\dbF}(\O;D([0,T];\dbR^n)) \t
L^{2}_{\dbF}(\O;L^2(0,T;\dbR^n))$ of equation (\ref{sm1}) is the
unique strong solution of this equation, and $y(\cdot) \in
L^{2}_{\dbF}(\O;C([0,T];\dbR^n))$.

\ms

ii) If $\ds
y_T-\int_0^Tf(s)ds-\dbE\left(y_T-\int_0^Tf(s)ds\right)\in \left(
L_{\dbW,\cF_T}^2(\O;\dbR^n)\right)^\perp$, then the transposition
solution $(y(\cdot), Y(\cdot)) \in L^{2}_{\dbF}(\O;D([0,T];\dbR^n))
\t L^{2}_{\dbF}(\O;L^2(0,T;\dbR^n))$ of equation (\ref{sm1}) is
given by the following
 \bel{9e6}
 \left\{
 \ba{ll}
 \ds y(t)=\dbE\left.\left(y_T-\int_t^Tf(s)ds\;\right|\;\cF_t\right),\\
 \ns\ds
 Y(\cd)=0.
 \ea\right.
 \ee
\ep

{\bf Proof.} i) By definition of $L_{\dbW,\cF_T}^2(\O;\dbR^n)$, one
can find a $\cl{Y}(\cd)\in L^2_{\dbF}(\O;L^2(0,T;\dbR^n))$ such that
 $$
  y_T-\int_0^Tf(s)ds-\dbE\left(y_T-\int_0^Tf(s)ds\right)=\int_0^T\cl{Y}(s)dw(s).
 $$
Then, put
 $$
 \bar y(t)=\dbE\left(y_T-\int_0^Tf(s)ds\right)+\int_0^tf(s)ds+\int_0^t\cl{Y}(s)dw(s).
 $$
It is clear that $(\bar y(\cd),\cl{Y}(\cd))\in
L^{2}_{\dbF}(\O;C([0,T];\dbR^n)) \t L^2_{\dbF}(\O;L^2(0,T;\dbR^n))$
is a strong solution for the linear BSDE (\ref{sm1}). Clearly,
$(\bar y(\cd),\cl{Y}(\cd))$ is also a transposition solution of
(\ref{sm1}). Hence, by the uniqueness of the transposition solution
for (\ref{sm1}), it follows that $(\bar y(\cd),\cl{Y}(\cd))=(
y(\cd),Y(\cd))$.

\ms

ii) Choosing $t=0$, $\eta=0$, $u(\cd)=0$ and $v(\cd)=Y(\cd)$ in
equation (\ref{system2}), we get $\ds z(\tau)=\int_0^\tau
Y(s)dw(s)$. Hence, by definition, (\ref{sol1}) is now specialized as
 \bel{9e5}
\mathbb{E}\Big\langle \int_0^T Y(t)dw(t),y_T\Big\rangle   =
\mathbb{E}\int_0^T \Big\langle \int_0^t Y(s)dw(s),f(t)\Big\rangle dt
+ \mathbb{E} \int_0^T|Y(t)|^2 dt.
 \ee
Noting that
 $$
 \ba{ll}\ds
 \mathbb{E}\int_0^T \Big\langle \int_0^t Y(s)dw(s),f(t)\Big\rangle
 dt\\\ns\ds=\mathbb{E}\int_0^T \Big\langle \int_0^T Y(s)dw(s),f(t)\Big\rangle
 dt-\mathbb{E}\int_0^T \Big\langle \int_t^T
 Y(s)dw(s),f(t)\Big\rangle\\\ns\ds
 =\mathbb{E}\Big\langle \int_0^T Y(t)dw(t),\int_0^Tf(s)ds\Big\rangle
 \ea
 $$
and that
 $$
\mathbb{E}\Big\langle \int_0^T
Y(t)dw(t),\dbE\left(y_T-\int_0^Tf(s)ds\right)\Big\rangle=0,
 $$
we conclude from (\ref{9e5}) that
  \bel{9e6-1}
\mathbb{E}\Big\langle \int_0^T Y(t)dw(t),\
y_T-\int_0^Tf(s)ds-\dbE\left(y_T-\int_0^Tf(s)ds\right)\Big\rangle =
\mathbb{E} \int_0^T|Y(t)|^2 dt.
 \ee
Now by our assumption that $\ds
y_T-\int_0^Tf(s)ds-\dbE\left(y_T-\int_0^Tf(s)ds\right)\in \left(
L_{\dbW,\cF_T}^2(\O;\dbR^n)\right)^\perp$, it follows from
(\ref{9e6-1}) that $Y(\cd)=0$.

Next, choosing $t=0$, $\eta=0$ and $v(\cd)=0$, and $u(\cd)\in
L^2_{\dbF}(\O;L^1(0,T;\dbR^n))$ (arbitrarily) in equation
(\ref{system2}), we get $\ds z(\tau)=\int_0^\tau u(s)ds$. Hence, by
definition, (\ref{sol1}) is now specialized as
 $$
\mathbb{E}\Big\langle \int_0^T u(t)dt,y_T\Big\rangle   =
\mathbb{E}\int_0^T \Big\langle \int_0^t u(s)ds,f(t)\Big\rangle dt +
\mathbb{E}  \int_0^T\langle u(t),y(t)\rangle dt.
 $$
Hence,
$$
\mathbb{E}\int_0^T\Big\langle  u(t),\ y(t)-y_T+\int_t^T
u(s)ds\Big\rangle dt  = 0.
 $$
This gives
 \bel{9e7}
\mathbb{E}\int_0^T\Big\langle  u(t),\
y(t)-\dbE\left.\left(y_T-\int_t^Tf(s)ds\;\right|\;\cF_t\right)\Big\rangle
dt  = 0,\q\forall\; u(\cd)\in L^2_{\dbF}(\O;L^1(0,T;\dbR^n)).
 \ee
Now, the first equality in (\ref{9e6}) follows from (\ref{9e7}).
This completes the proof of Proposition \ref{3p1}.
\endpf

\br\label{3r3}
Proposition \ref{3p1} ii) justifies our transposition solution.
Indeed, when $ L_{\dbW,\cF_T}^2(\O;\dbR^n)$ is a proper subspace of
$ L_{0,\cF_T}^2(\O;\dbR^n)$ and $\ds
y_T-\int_0^Tf(s)ds-\dbE\left(y_T-\int_0^Tf(s)ds\right)\in \left(
L_{\dbW,\cF_T}^2(\O;\dbR^n)\right)^\perp$, it is easy to show that
the transposition solution (\ref{9e6}) is NOT a strong solution of
equation (\ref{sm1}).
\er

\br
As far as we know, there exists no any satisfactory characterization
on $ L_{\dbW,\cF_T}^2(\O;\dbR^n)$. Especially, it seems to us that
it is not very clear when $ L_{\dbW,\cF_T}^2(\O;\dbR^n)$ is a proper
subspace of $ L_{0,\cF_T}^2(\O;\dbR^n)$. Of course, it is easy to
see that $ L_{\dbW,\cF_T}^2(\O;\dbR^n)= L_{0,\cF_T}^2(\O;\dbR^n)$
implies that the Martingale Representation Theorem holds.
\er

\section{Well-posedness of semilinear BSDEs}

The purpose of this section is to establish the following
well-posedness result for the semilinear BSDE (\ref{system1}).

\begin{theorem}\label{theorem1}
For any given $y_T \in L^2_{\cF_T}(\O)$, equation (\ref{system1})
admits a unique transposition solution $(y(\cdot), Y(\cdot))\in
L^{2}_{\dbF}(\O;D([0,T];\dbR^n)) \t
L^{2}_{\dbF}(\O;L^2(0,T;\dbR^n))$. Furthermore, there is a constant
$C>0$, depending only on $K$ and $T$, such that
 \bel{et11}
 \ba{ll}\ds
 |(y(\cdot),   Y(\cdot))|_{
L^{2}_{\dbF}(\O;D([0,T];\dbR^n)) \t
L^{2}_{\dbF}(\O;L^2(0,T;\dbR^n))}\\\ns \ds\le C\left[|f(\cd,0,0)|_{
L^2_{\dbF}(\O;L^1(0,T;\dbR^n))}+|y_T|_{
L^2_{\cF_T}(\O;\dbR^n)}\right].
 \ea
 \ee
\end{theorem}

{\bf Proof.} Fix any $T_1\in [0,T]$. For any $(p(\cdot), P(\cdot))
\in L^{2}_{\dbF}(\O;D([T_1,T];\dbR^n)) \t
L^{2}_{\dbF}(\O;L^2(T_1,T;\dbR^n))$, we consider the following
equation:
\begin{eqnarray}\label{iterative2}
\left\{
\begin{array}{lll}
\ds dy = f(t,p(t),P(t)) dt + Ydw & \mbox{ in } [T_1,T],\\
\ns\ds y(T)=y_T.
\end{array}
\right.
\end{eqnarray}
By condition (\ref{Lm1}) and Theorem \ref{the1}, equation
(\ref{iterative2}) admits a transposition solution $(y(\cdot),
Y(\cdot)) \in L^{2}_{\dbF}(\O;$ $D([T_1,T];\dbR^n)) \t
L^{2}_{\dbF}(\O;L^2(T_1,T;\dbR^n))$. This defines a map $F$ from
$L^{2}_{\dbF}(\O;D([T_1,T];\dbR^n)) \t
 L^{2}_{\dbF}(\O;$ $L^2(T_1,T;\dbR^n))$ into itself by
$F(p(\cdot), P(\cdot))  = (y(\cdot), Y(\cdot))$.

We claim that the map $F$ is contractive provided that $T-T_1$ is
small enough. Indeed, for another $(\hat p(\cdot), \hat P(\cdot))
\in L^{2}_{\dbF}(\O;D([T_1,T];\dbR^n)) \t
L^{2}_{\dbF}(\O;L^2(T_1,T;\dbR^n))$, we define $(\hat y(\cdot), \hat
Y(\cdot))=F(\hat p(\cdot), \hat P(\cdot))$. Put
 $$
 \tilde y(\cd)=y(\cd)-\hat y(\cd),\q \widetilde Y(\cd)=Y(\cd)-\hat
 Y(\cd),\q \tilde f(\cd)=f(\cd,p(\cd),P(\cd))-f(\cd,\hat p(\cd),\hat
 P(\cd)).
 $$
Clearly, $(\tilde y(\cdot), \widetilde Y(\cdot))$ solves the
following equation
 \begin{eqnarray}\label{iterative3}
\left\{
\begin{array}{lll}
\ds d \tilde y = \tilde f(t) dt +\widetilde Ydw & \mbox{ in } [T_1,T],\\
\ns\ds \tilde y(T)=0.
\end{array}
\right.
\end{eqnarray}
By condition (\ref{Lm1}), it is easy to see that $\tilde f(\cd)\in
L^2_{\dbF}(\O;L^1(T_1,T;\dbR^n))$ and
 \bel{et3}
 \ba{ll}\ds
 |\tilde f(\cd)|_{L^2_{\dbF}(\O;L^1(T_1,T;\dbR^n))}\\
  \ns\ds\le K\left[|p(\cd)-\hat
 p(\cd)|_{L^2_{\dbF}(\O;L^2(T_1,T;\dbR^n))}+|P(\cd)-\hat
 P(\cd)|_{L^2_{\dbF}(\O;L^1(T_1,T;\dbR^n))}\right]\\
 \ns
 \ds\le K\big(T-T_1+\sqrt{T-T_1}\big)\left[|p(\cd)-\hat
 p(\cd)|_{L^{2}_{\dbF}(\O;D([T_1,T];\dbR^n)) }+|P(\cd)-\hat
 P(\cd)|_{L^2_{\dbF}(\O;L^2(T_1,T;\dbR^n))}\right].
 \ea
 \ee
Applying Theorem \ref{the1} to equation (\ref{iterative3}) and
noting (\ref{et3}), it follows that there is a constant $C$,
depending only on $T$, such that
  \bel{et8}
 \ba{ll}\ds
 |(\tilde y(\cdot), \tilde  Y(\cdot))|_{
L^{2}_{\dbF}(\O;D([T_1,T];\dbR^n)) \t
L^{2}_{\dbF}(\O;L^2(T_1,T;\dbR^n))}\le
C |\tilde f(\cd)|_{ L^2_{\dbF}(\O;L^2(T_1,T;\dbR^n))}\\
\ns \ds\le CK\big(T-T_1+\sqrt{T-T_1}\big)\left[|p(\cd)-\hat
 p(\cd)|_{L^{2}_{\dbF}(\O;D([T_1,T];\dbR^n)) }+|P(\cd)-\hat
 P(\cd)|_{L^2_{\dbF}(\O;L^2(T_1,T;\dbR^n))}\right].
 \ea
 \ee
One may choose $T_1$ so that $CK\big(T-T_1+\sqrt{T-T_1}\big)<1$, and
hence $F$ is a contractive map.

By the Banach fixed point theorem, $F$ has a fixed point $(y(\cdot),
Y(\cdot))\in L^{2}_{\dbF}(\O;D([T_1,T];\dbR^n)) \t
L^{2}_{\dbF}(\O;L^2(T_1,T;\dbR^n))$. It is clear that $(y(\cdot),
Y(\cdot))$ is a transposition solution to the following equation:
 \begin{eqnarray}\label{iterative5}
\left\{
\begin{array}{lll}
\ds dy = f(t,y(t),Y(t)) dt + Ydw & \mbox{ in } [T_1,T],\\
\ns\ds y(T)=y_T.
\end{array}
\right.
\end{eqnarray}
Using again condition (\ref{Lm1}) and similar to (\ref{et3}), we see
that $ f(\cd,y(\cd),Y(\cd))\in L^2_{\dbF}(\O;L^1(T_1,T;\dbR^n))$ and
  \bel{et008}
 \ba{ll}\ds
  |f(\cd,y(\cd),Y(\cd))|_{ L^2_{\dbF}(\O;L^1(T_1,T;\dbR^n))}\\
  \ns\ds\le |f(\cd,0,0)|_{ L^2_{\dbF}(\O;L^1(T_1,T;\dbR^n))}+K\left[|y(\cd)|_{L^2_{\dbF}(\O;L^1(T_1,T;\dbR^n))}
  +|Y(\cd)|_{L^2_{\dbF}(\O;L^1(T_1,T;\dbR^n))}\right]\\
\ns \ds\le |f(\cd,0,0)|_{
L^2_{\dbF}(\O;L^1(T_1,T;\dbR^n))}\!+\!K\big(T\!\!-\!\!T_1\!+\!\sqrt{T\!-\!T_1}\big)\left[|y(\cd)|_{L^{2}_{\dbF}(\O;D([T_1,T];\dbR^n))
}\!+\!|Y(\cd)|_{L^2_{\dbF}(\O;L^2(T_1,T;\dbR^n))}\right].
 \ea
 \ee
Applying Theorem \ref{the1} to equation (\ref{iterative5}) and
noting (\ref{et008}), we find that
  \bel{et9}
 \ba{ll}\ds
  |(y(\cdot),   Y(\cdot))|_{
L^{2}_{\dbF}(\O;D([T_1,T];\dbR^n)) \t L^{2}_{\dbF}(\O;L^2(T_1,T;\dbR^n))}\\
\ns \ds\le C\left[
 |f(\cd,y(\cd),Y(\cd))|_{ L^2_{\dbF}(\O;L^1(T_1,T;\dbR^n))}+|y_T|_{
L^2_{\cF_T}(\O;\dbR^n)}\right]\\
\ns \ds\le C\left[K\big(T-T_1+\sqrt{T-T_1}\big)|(y(\cdot),
Y(\cdot))|_{ L^{2}_{\dbF}(\O;D([T_1,T];\dbR^n)) \t
L^{2}_{\dbF}(\O;L^2(T_1,T;\dbR^n))}\right.\\
\q\ds\left.+|f(\cd,0,0)|_{ L^2_{\dbF}(\O;L^1(T_1,T;\dbR^n))}+|y_T|_{
L^2_{\cF_T}(\O;\dbR^n)}\right].
 \ea
 \ee
Noting that $K\big(T-T_1+\sqrt{T-T_1}\big)<1$, by (\ref{et9}), we
get
   \bel{et10}
 \ba{ll}\ds
 |(y(\cdot),   Y(\cdot))|_{
L^{2}_{\dbF}(\O;D([T_1,T];\dbR^n)) \t
L^{2}_{\dbF}(\O;L^2(T_1,T;\dbR^n))}\le C\left[|f(\cd,0,0)|_{
L^2_{\dbF}(\O;L^1(0,T;\dbR^n))}+|y_T|_{
L^2_{\cF_T}(\O;\dbR^n)}\right].
 \ea
 \ee

Repeating the above argument step by step, we obtain the
transposition solution of equation (\ref{system1}) on $[0,T]$. The
uniqueness of this solution is obvious. The desired estimate
(\ref{et11}) follows from (\ref{et10}). This completes the proof of
Theorem \ref{theorem1}.\endpf

\br\label{4r1}
For the transposition solution $(y(\cd),Y(\cd))$ to equation
(\ref{system1}), put
 \bel{9e21}
 M(t)=y(t)-\int_0^tf(s,y(s),Y(s))ds,\qq\forall\; t\in [0,T].
  \ee
Thanks to Remark \ref{3r1}, it is easy to see that $M(\cd)$ is a
$\cF_t$-martingale, and $(y(\cd),M(\cd),Y(\cd))$ satisfies the
following equation
  \bel{9e3}
 y (t) = y_T-\int_t^Tf(s,y(s),Y(s))ds +M(t)-M(T),\qq\forall\;t\in
 [0,T].
 \ee
Equation (\ref{9e3}) can be regarded as a corrected form of equation
(\ref{system1}).
\er

Stimulating by Remark \ref{4r1}, we introduce the following notion
for solution of equation (\ref{system1}).

\begin{definition}\label{dn}
We call $(y(\cd),M(\cd),Y(\cd)) \in \Upsilon \t
L^2_{\dbF}(\O;L^2(0,T;\dbR^n))$ to be a corrected solution of
equation (\ref{system1}) if  $(y(\cdot), Y(\cdot))$ is a
transposition solution of this equation, and (\ref{9e3}) holds.
\end{definition}

As a consequence of Theorem \ref{theorem1} and Remark \ref{4r1}, it
is easy to prove the following result.

\begin{corollary}\label{th1}
For any given $y_T \in L^2_{\cF_T}(\O)$, equation (\ref{system1})
admits a unique corrected solution $(y(\cd),M(\cd),Y(\cd)) \in
\Upsilon \t L^2_{\dbF}(\O;L^2(0,T;\dbR^n))$. Furthermore, there is a
constant $C>0$, depending only on $K$ and $T$, such that
  $$
 |(y(\cd),M(\cd),Y(\cd))|_{\Upsilon \t L^2_{\dbF}(\O;L^2(0,T;\dbR^n))}\le C\left[|f(\cd,0,0)|_{
L^2_{\dbF}(\O;L^1(0,T;\dbR^n))}+|y_T|_{
L^2_{\cF_T}(\O;\dbR^n)}\right].
  $$
\end{corollary}

\br
Clearly, for the corrected solution $(y(\cd),M(\cd),Y(\cd)) \in
\Upsilon \t L^2_{\dbF}(\O;L^2(0,T;\dbR^n))$ of equation
(\ref{system1}) obtained in Corollary \ref{th1}, the first two
components satisfy (\ref{9e21}). Furthermore, if the filtration
$\dbF$ is the natural one, then the last two components of this
solution satisfies (\ref{7e5}).
\er

\br
Using the method developed in \cite{KH}, one can find a unique
solution
 $
(\tilde y(\cd), \wt Q(\cd),\wt{Y}(\cd))\in
L^{2}_{\dbF}(\O;D([0,T];\dbR^n)) \t
\left(\cM_{0,\dbM,\dbF}^2([0,T];\dbR^n)\right)^\perp\times
L^2_{\dbF}(\O;L^2(0,T;\dbR^n)) $ satisfying the following equation
\begin{equation}\label{third3}
\tilde y (t) = y_T+\wt Q(t)-\wt Q(T)-\int_t^Tf(s,\tilde y(s),\wt
Y(s))ds -\int_t^T \wt Y(s) dw
 (s),\qq\forall\;t\in [0,T].
 \end{equation}
Equation (\ref{third3}) can be regarded as another corrected form of
equation (\ref{system1}). Using It\^o's formula to $\lan z(t),\tilde
y(t)\ran$ (Recall (\ref{system2}) for $z(\cd)$), and noting the
strong orthogonality of the martingales $\wt Q(\cd)$ and
$\ds\int_0^\cdot \wt Y(s)dw(s)$ (which follows from some deep
results in martingale theory, e.g., \cite[Chapter VIII]{DM}), one
can show that
 \bel{sol1-2}
 \ba{ll}
 \ds\mathbb{E}\lan z(T),y_T\ran - \mathbb{E}\lan
\eta,\tilde y(t)\ran\\
 \ns
\ds = \mathbb{E}\int_t^T \lan z(\tau),f(\tau,\tilde
 y(\tau),\wt Y(\tau))\ran d\tau + \mathbb{E}\int_t^T \lan
u(\tau),\tilde y(\tau)\ran d\tau + \mathbb{E} \int_t^T\lan v(\tau),
\wt Y(\tau)\ran d\tau
 \ea
\ee
holds for all $t\in [0,T]$, $u(\cdot)\in
L^2_{\dbF}(\O;L^1(t,T;\dbR^n))$, $v(\cdot)\in
 L^2_{\dbF}(\O;L^2(t,T;\dbR^n))$ and $\eta\in
 L^2_{\cF_t}(\O;\dbR^n)$. Hence, $
(\tilde y(\cd), \wt{Y}(\cd))$ is also a transposition solution of
(\ref{system1}). By the uniqueness result in Corollary \ref{th1}, we
conclude that
 $$(y(\cd),M(\cd),Y(\cd)) = (\tilde y(\cd),M(0)+\int_0^\cdot \wt Y(s)dw(s)+\wt Q(\cd),\wt Y(\cd)).
 $$
Nevertheless, as we explained before, in some sense, our method
seems to be more flexible for the general filtration than the
existing ones.
\er

 \br
In some sense, the transposition/corrected solution for BSDEs is in
spirit close to the distribution solution for partial differential
equations. It is then very natural to study the further regularity
for the transposition solution $(y(\cd),Y(\cd))$ (or corrected
solution $(y(\cd),M(\cd),Y(\cd))$) of equation (\ref{system1}).
\er

\section{Comparison theorem for transposition solutions}

In this section, we show a comparison theorem for transposition
solutions of the semilinear BSDE (\ref{system1}) in one dimension,
i.e., $n=1$.

We will go a little further. Besides equation (\ref{system1}) (with
$n=1$), we consider also the following BSDE:
\begin{eqnarray}\label{c system2}
\left\{
\begin{array}{lll}
\ds d\bar{y} (t) = \bar{f}(t,\bar{y}(t),\cl{Y}(t))dt + \cl{Y}(t) dw (t) & \mbox{ in } [0,T], \\
 \ns\ds \bar{y}(T) = \bar{y}_T.
\end{array}
\right.
 \end{eqnarray}
Here $\bar{y}_T \in L^2_{\cF_T}(\O;\dbR)$, $\bar{f}(\cd,\cd,\cd)$ is
supposed to satisfy $\bar{f}(\cd,0,0)\in
L^2_{\dbF}(\O;L^1(0,T;\dbR))$ and,
  \bel{cLm2}
 |\bar{f}(t,p_1,q_1)-\bar{f}(t,p_2,q_2)|\le
K(|p_1-p_2|+|q_1-q_2|),\q t\in [0,T]\hb{ a.s.},
\forall\;p_1,p_2,q_1,q_2\in \dbR.
 \ee
By Theorem \ref{theorem1}, equation (\ref{c system2}) admits a
unique transposition solution $(\bar y(\cdot), \cl{Y}(\cdot))\in
L^{2}_{\dbF}(\O;D([0,T];$ $\dbR)) \t
L^{2}_{\dbF}(\O;L^2(0,T;\dbR))$.

We have the following result.

\begin{theorem}\label{comparison th}
If $y_T \geq \bar{y}_T$ a.s., and for any $a,b\in\dbR$, and a.e.
$t\in [0,T]$,
\begin{equation}\label{c condition1}
f(t,a,b) \leq \bar{f}(t,a,b)\ \mbox{ a.s.,}
\end{equation}
then, for any $t\in [0,T]$,
\begin{equation}\label{c result1}
y (t)\geq \bar{y}(t)\ \mbox{ a.s.}
\end{equation}
Moreover, $y(t) = \bar{y}(t)$ a.s., for some $t\in [0,T]$ if and
only if $y_T = \bar{y}_T$ a.s., and that $f(s,\bar{y}(s),\cl{Y}(s))
= \bar{f}(s,\bar{y}(s),\cl{Y}(s))$ a.s. for a.e. $s\in [t,T]$.
\end{theorem}

{\bf Proof}. The idea of our proof is very close to that of
\cite[Theorem 2.2]{KPQ}. Put $\hat{y} = y - \bar{y}$ and $\hat Y = Y
- \cl{ Y}$. It is clear that $(\hat y,\hat Y)$ is a transposition
solution of the following equation
\begin{eqnarray}\label{c system3}
\left\{
\begin{array}{lll}
\ds d\hat y (t) = (a(t)\hat y(t) + b(t)\hat Y(t) + h(t))dt + \hat Y(t) dw(t) & \mbox{ in } [0,T], \\
 \ns\ds \hat y(T) = y_T - \bar y_T,
\end{array}
\right.
 \end{eqnarray}
where
\begin{eqnarray*}
a(t) = \left\{
\begin{array}{lll}
\ds\frac{f(t,y(t),Y(t)) - f(t,\bar y(t),Y(t))}{y(t)-\bar y(t)}, &
y(t) \neq \bar y(t),\\
\ns\ds 0, & y(t) = \bar y(t),
\end{array}
\right.
\end{eqnarray*}
\begin{eqnarray*}
b(t) = \left\{
\begin{array}{lll}
\ds\frac{f(t,\bar y(t),Y(t)) - f(t,\bar y(t),\cl{ Y}(t))}{\ds
Y(t)-\cl{ Y}(t)}, &
Y(t) \neq\cl{ Y}(t),\\
\ns\ds 0, & Y(t) = \cl{ Y}(t),
\end{array}
\right.
\end{eqnarray*}
and
$$
h(t) = f(t,\bar y(t),\cl{ Y}(t)) - \bar f(t,\bar y(t),\cl{
Y}(t))\leq 0.
$$
From (\ref{Lm1}) and (\ref{cLm2}), we see that $|a(t)|\leq K$ and
$|b(t)|\leq K$ a.s., for a.e. $t\in [0,T]$.

Now, for any $t\in[0,T]$, we consider the following (forward)
stochastic differential equation
\begin{eqnarray}\label{cdsystem3}
\left\{
\begin{array}{lll}
\ds d q (s) = -a(s)q(s) dt - b(s)q(s) dw(s) & \mbox{ in } [t,T], \\
 \ns\ds q(t) = \varsigma,
\end{array}
\right.
 \end{eqnarray}
where $\varsigma \in L^2_{\cF_t}(\O;\dbR)$ satisfying $\varsigma
\geq 0$ a.s. It is easy to see that
\begin{eqnarray}\label{c eq1}q(s) = \varsigma
\exp\Big\{-\int_t^s a(\tau)d\tau - \frac{1}{2}\int_t^s b(\tau)d\tau
- \int_t^s b(\tau)dw(\tau) \Big\}\geq 0.
\end{eqnarray}
 Since $(\hat
y,\hat Y)$ is the transposition solution of equation (\ref{c
system3}), by Definition \ref{def of solution}, it follows that
\begin{eqnarray*}
&\,&\dbE\big( \hat y(T)q(T) \big)- \dbE \big(\hat y(t)
\varsigma\big) = \dbE\int_t^T [a(s)\hat y(s)
+ b(s)\hat Y(s) + h(s)]q(s)ds \nonumber\\
&&\,\qq\qq\qq\qq\qq\qq - \dbE \int_t^T \hat y(s) a(s)q(s)ds - \dbE
\int_t^T \hat Y(s) b(s)q(s)ds,
\end{eqnarray*}
from which we clonclude that
\begin{eqnarray}\label{c eq2}
\dbE\big( \hat y(t) \varsigma \big)= \dbE \big(\hat y(T)q(T)\big) -
\dbE\int_t^T h(s)q(s)ds \geq 0,
\end{eqnarray}
for any $\varsigma \in L^2_{\cF_t}(\O;\dbR)$ such that $\varsigma
\geq 0$ a.s. Therefore, we see that $\hat y (t) \geq 0$ a.s., which
means that $y(t)\geq \bar y(t)$ a.s.

Choosing $\varsigma = 1$ in (\ref{cdsystem3}), from (\ref{c eq1}),
it is easy to see that $q(s)>0$ for any $s\in [t,T]$. By (\ref{c
eq2}), we obtain that

\begin{eqnarray*}
\dbE \hat y(t) = \dbE \big(\hat y(T)q(T)\big) - \dbE\int_t^T
h(s)q(s)ds \geq 0.
\end{eqnarray*}
If $\hat y(t)=0$ a.s., it follows that $\dbE \big(\hat
y(T)q(T)\big)=0$ and $\ds \dbE\int_t^T  h(s)q(s)ds = 0$. Since
$q(s)>0$ for any $s\in [t,T]$, we have that $\hat y(T)=0$ a.s. and
$h(s) = 0$ a.s. for a.e. $s\in [t,T]$, which leads to $y_T =
\bar{y}_T$ a.s., and that $f(s,\bar{y}(s),\cl{Y}(s)) =
\bar{f}(s,\bar{y}(s),\cl{Y}(s))$ a.s., for a.e. $s\in [t,T]$. \endpf

\section*{Acknowledgement}

This work is supported by the NSFC under grants 10831007 and
60974035,  and the project MTM2008-03541 of the Spanish Ministry of
Science and Innovation. The second author acknowledges gratefully
Professors Zhenqing~Chen, Shige~Peng, Jia-An Yan, Jiongmin~Yong and
Xunyu~Zhou for stimulating discussions, and Dr. Mingyu~Xu for
pointing out reference \cite{LLQ}.




\end{document}